%% file: PIB_ver1.tex
\documentclass[12pt,final]{article}

\input paper.tex

\title{The characterization of cyclic cubic fields with power integral bases}
\author{Tomokazu Kashio\thanks{Tokyo University of Science, \protect\url{kashio_tomokazu@ma.noda.tus.ac.jp}} \\
Ryutaro Sekigawa\thanks{Tokyo University of Science, \protect\url{sekigawa.r@gmail.com}}}

\begin{document}

\maketitle

\begin{abstract}
We provide an equivalent condition for the monogenity of the ring of integers of any cyclic cubic field.
Although a large part of the main results is covered by the classical one of Gras, we write the condition more explicitly. 
First, we show that if a cyclic cubic field is monogenic then it is a simplest cubic field $K_t$ defined by Shanks' cubic polynomial $f_t(x):=x^3-tx^2-(t + 3)x-1$ with $t \in \mathbb Z$. 
Then we give an equivalent condition for when $K_t$ is monogenic, which is explicitly written in terms of $t$. 
\end{abstract}

\section{Main results}

We say a number field $K$ has a {\it power integral basis} if its ring of integers $\mathcal O_K$ is generated by one element: $\mathcal O_K=\mathbb Z[\gamma]$.
We say also that $K$ or $\mathcal O_K$ is {\it monogenic}.
Let $K$ be an abelian field of prime degree $p$. 
It is known that $K$ has a power in integral basis when $p=2$, 
and does not have when $p\geqq 5$ unless $K= \mathbb Q(\zeta_{2p+1}+\zeta_{2p+1}^{-1})$ and $2p+1$ is a prime \cite{Gr4}.
In this sense, the case of $p=3$ (i.e., cyclic cubic fields) is of great interest.
In fact, Dummit-Kisilevsky \cite{DK} showed that there exist infinitely many cyclic cubic fields with or without power integral bases, respectively.
Archinard \cite{Ar} and Gras \cite{Gr1,Gr2,Gr3} also studied this case. 

Now we introduce the Shanks cubic polynomials \cite{Sh}
\begin{align*}
f_t(x)=x^3-tx^2-(t + 3)x-1 \quad (-1\leqq t \in \mathbb Z). 
\end{align*}
Let $\theta_t$ be a root of $f_t(x)=0$ and put $K_t:=\mathbb Q(\theta_t)$. We assume that $-1\leqq t$ since we have $K_t=K_{-(t+3)}$. The followings are well-known facts:
\begin{itemize}
\item $K_t/\mathbb Q$ is a cyclic cubic extension. We put 
\begin{align*}
G:=\mathrm{Gal}(K_t/\mathbb Q)=\langle \sigma \rangle=\{\mathrm{id},\sigma,\sigma^2\}.
\end{align*}
\item $\theta_t$ is a unit $\in \mathcal O_{K_t}^\times$ satisfying that (if necessary by replacing $\sigma$ with $\sigma^2$)
\begin{align*}
\sigma(\theta_t)=-\frac{1+\theta_t}{\theta_t},&& \sigma^2(\theta_t)=\frac{-1}{1+\theta_t},  
&&1+\theta_t+\theta_t\sigma(\theta_t)=0,  &&\mathit{N}(\theta_t)=1,  \quad \mathit{Tr}(\theta_t)=t.
\end{align*}
\item Let $c_K$ and $d_K$ denote the conductor and the discriminant, respectively. We have  
\begin{align*}
&c_{K_t}=\sqrt{d_{K_t}} \mid \Delta_t:=t^2 + 3t + 9.
\end{align*}
\end{itemize}
The relation $c_K=\sqrt{d_K}$ follows from the conductor-discriminant formula.
We abbreviate as $\mathit{N}:=\mathit{N}_{K/\mathbb Q}$, $\mathit{Tr}:=\mathit{Tr}_{K/\mathbb Q}$.
We call a cyclic cubic field of the form $K=K_t$ ($-1\leqq t \in \mathbb Z$) a {\it simplest cubic field}.

We state the main result in this paper.
Let $v_p(n)$ denote the order at a prime $p$ of an integer $n$ defined by 
\begin{align*}
n=p^{v_p(n)}n_0, \quad (p,n_0)=1.
\end{align*}

\begin{thm} \label{thm1}
Let $K$ be a cyclic cubic field. 
The followings are equivalent.
\begin{enumerate}
\item $K$ has a power integral basis.
\item There exists $t$ satisfying that $K=K_t$ and 
\begin{align*}
\frac{\Delta_t}{c_K} \in \mathbb N^3:=\left\{n^3 \mid n \in \mathbb N\right\}.
\end{align*}
\item There exists $t$ satisfying that $K=K_t$ and 
\begin{align*}
t \not \equiv 3,21 \bmod 27, && v_p(\Delta_t) \not\equiv 2 \bmod 3 \text{ for all }p \neq 3.
\end{align*}
\end{enumerate}
In such a case, a power integral basis $\gamma \in \mathcal O_K$ is given by
\begin{align} \label{pib}
\gamma:=
\frac{\theta_t-a}{\sqrt[3]{\frac{\Delta_t}{c_K}}} \text{ for } a \in \mathbb Z \text{ with } a\equiv \frac{t}{3} \bmod \sqrt[3]{\frac{\Delta_t}{c_K}}.
\end{align}
Here the condition $a \equiv \frac{t}{3} \bmod \sqrt[3]{\frac{\Delta_t}{c_K}}$ means the following.
We have $3 \mid t$ when $3 \mid \sqrt[3]{\frac{\Delta_t}{c_K}}$ by {\rm\cref{3}-(i)} in {\rm\S \ref{ramify}}. Hence, strictly speaking, we take an integer $a$ satisfying     
\begin{align*}
\begin{cases}
a\equiv \dfrac{t}{3} \bmod \sqrt[3]{\frac{\Delta_t}{c_K}} & \left(3 \mid \sqrt[3]{\frac{\Delta_t}{c_K}}\right), \\[5pt]
3a\equiv t \bmod \sqrt[3]{\frac{\Delta_t}{c_K}} & \left(3 \nmid \sqrt[3]{\frac{\Delta_t}{c_K}}\right).
\end{cases}
\end{align*}
\end{thm}

\begin{rmk}
The condition in {\rm\cref{thm1}-(iii)} or {\rm\cref{crl}-(iii)} below is given in terms of $t$ (and $\Delta_t=t^2+3t+9$).
Namely, we can study the monogenity without studying the number field $K_t$ itself.
\end{rmk}

Let us explain the results of the paper more precisely.
First we obtain the following characterization (\cref{thm2}) of cyclic cubic fields having power integral bases.
We prepare some notations.
For a number field $k$, we denote by $I_k,P_k$ the group of all fractional ideals, all principal ideals, respectively.
Additionally, we put
\begin{align*}
I_K^G:=\left\{\mathfrak A \in I_K \mid \sigma(\mathfrak A)=\mathfrak A\right\}, &&
P_K^G:=I_K^G \cap P_K
\end{align*}
to be the group of all ambiguous ideals, all principal ambiguous ideals, respectively.
Let $\mathfrak C_K$ be a unique integral ideal satisfying 
\begin{align*}
\mathit{N} \mathfrak C_K=c_K,
\end{align*}
which is given by \cref{fortriv}-(i) in \S \ref{mvss}. 
We easily see that
\begin{itemize}
\item $\mathfrak C_K \in I_K^G$. In particular, if $\mathfrak C_K$ is principal, then $\mathfrak C_K \in P_K^G$.
\item Let $K=K_t$. Then $(\theta_t-\sigma(\theta_t)) \in P_K^G$.
\end{itemize}
In fact, by definition, we have $\sigma(\mathfrak C_K)=\mathfrak C_{\sigma(K)}=\mathfrak C_K$. 
The latter one follows from an explicit calculation:
\begin{align} \label{amb}
\frac{\sigma(\theta_t)-\sigma^2(\theta_t)}{\theta_t-\sigma(\theta_t)}
=\frac{-\frac{1+\theta_t}{\theta_t}-\frac{-1}{1+\theta_t}}{\theta_t+\frac{1+\theta_t}{\theta_t}}
=\frac{-1}{1+\theta_t}=\sigma^2(\theta_t) \in \mathcal O_K^\times.
\end{align}
Additionally we consider the natural projection
\begin{align*}
P_K^G \ra P_K^G/P_\mathbb Q,\ \mathfrak A \mt \overline{\mathfrak A}:=\mathfrak A \bmod P_\mathbb Q.
\end{align*}
We also consider the $1$st cohomology group 
\begin{align*}
H^1(G,\mathcal O_K^\times)=\left\{u\in \mathcal O_K^\times \mid \mathit{N}(u)=1 \right\}\big/\left\{u^{\sigma-1} \mid u \in \mathcal O_K^\times\right\}
\end{align*}
and denote the class of $u\in \mathcal O_K^\times$ with $\mathit{N} (u)=1$ by $[u] \in H^1(G,\mathcal O_K^\times)$.

\begin{thm} \label{thm2}
Let $K$ be a cyclic cubic field. 
\begin{enumerate}
\item If $K$ has a power integral basis, then we have
\begin{itemize}
\item $K$ is a simplest cubic field.
\item $\mathfrak C_K$ is principal.
\end{itemize}
\item Assume that $K$ is a simplest cubic field and $\mathfrak C_K$ is principal, say $\mathfrak C_K=(\beta)$. The followings are equivalent.
\begin{enumerate}
\item $K$ has a power integral basis.
\item There exists $t$ satisfying that $K=K_t$ and $[\beta^{\sigma-1}]=[\theta_t] \in H^1(G,\mathcal O_K^\times)$. 
\item There exists $t$ satisfying that $K=K_t$ and $\overline{\mathfrak C_K}=\overline{(\theta_t-\sigma(\theta_t))} \in P_K^G/P_\mathbb Q$. 
\item There exists $t$ satisfying that $K=K_t$ and $\frac{\Delta_t}{c_K} \in \mathbb N^3$. 
\end{enumerate}
Note that $\beta^{\sigma-1}=\frac{\sigma(\beta)}{\beta} \in \mathcal O_K^\times$ since $\mathfrak C_K$ is an ambiguous ideal
and that the cohomology class $[\beta^{\sigma-1}] \in H^1(G,\mathcal O_K^\times)$ does not depend on the choice of $\beta$ since 
$[(\beta\epsilon )^{\sigma-1}]=[\beta^{\sigma-1}\epsilon^{\sigma-1}]=[\beta^{\sigma-1}]$ for $\epsilon \in \mathcal O_K^\times$.
\end{enumerate}
\end{thm}

Next, we give certain equivalent conditions for when $\mathfrak C_{K_t}$ is principal.
We also give an explicit condition in \cref{expl} in \S \ref{cvst}.

\begin{thm} \label{thm3}
Let $K=K_t$. For simplicity, 
assume that $K\neq \mathbb Q(\zeta_9+\zeta_9^{-1})$. 
The followings are equivalent.
\begin{enumerate}
\item $\mathfrak C_K$ is principal.
\item $\overline{\mathfrak C_K}=\overline{(\theta_t-\sigma(\theta_t))}$ or $\overline{(\theta_t-\sigma(\theta_t))^2} \in P_K^G/P_\mathbb Q$. 
\item $\frac{\Delta_t}{c_K}$ or $\frac{\Delta_t^2}{c_K} \in \mathbb N^3$. 
\end{enumerate}
Note that the condition {\rm(i)} is independent of the choice of $t$ with $K=K_t$, so are {\rm (ii)}, {\rm (iii)}. 
We also note that the exception $\mathfrak C_{\mathbb Q(\zeta_9+\zeta_9^{-1})}$ is a principal ideal ($=(\zeta_9+\zeta_9^{-1}+1)^2$). 
\end{thm}

Then, in \S \ref{pf}, we derive \cref{thm1} from Theorems \ref{thm2}, \ref{thm3} (and \cref{expl}), and some properties of $c_K,\Delta_t$ introduced in \S \ref{ramify}.

\begin{rmk}
The following characterization of $c_{K_t}$ summarizes the results in {\rm\S \ref{ramify}}:
\begin{align} \label{factor}
c_{K_t}=
\begin{cases}
\displaystyle \prod_{v_p(\Delta_t)\not \equiv 0 \bmod 3} p & (3 \nmid t \text{ or } t\equiv 12 \bmod 27), \\[3pt]
\displaystyle 3^2\prod_{p\neq 3,v_p(\Delta_t)\not \equiv 0 \bmod 3} p & (\text{otherwise}).
\end{cases}
\end{align}
We derive this expression in {\rm\S \ref{ramify}} although it seems to be well-known.
\end{rmk}

By \cite{Ok} or \cite[Theorem1.4]{Ho}, we have $K_t\neq K_{t'}$ if $t\neq t'$ except for  
\begin{align*}
K_{-1}=K_{5}=K_{12}=K_{1259}, && K_{0}=K_{3}=K_{54}(=\mathbb Q(\zeta_9+\zeta_9^{-1})), && K_{1}=K_{66}, && K_{2}=K_{2389}. 
\end{align*}
We see that these exceptions have power integral bases by applying \cref{thm1} for $t=-1,0,1,2$ (\cref{-1to70} in \S \ref{exm}).
Assume that $t\neq -1,0,1,2,3,5,12,54,66,1259,2389$. 
Then, in \cref{thm1}, \cref{thm2}, we may replace ``There exists $t$ satisfying that $K=K_t$ and $\sim$'' 
with ``For a unique $t$ satisfying $K=K_t$, we have $\sim$''. 
In particular, we have the following.

\begin{crl} \label{crl}
If a cyclic cubic field $K$ have a power integral basis, then it is a simplest cubic field, that is, there exists $t$ satisfying that $K=K_t$.
Moreover, the following are equivalent.
\begin{enumerate}
\item $K_t$ has a power integral basis.
\item $t \in \{ -1,0,1,2,3,5,12,54,66,1259,2389\}$ or $t$ satisfies that $\frac{\Delta_t}{c_{K_t}} \in \mathbb N^3$.
\item $t \in \{ -1,0,1,2,3,5,12,54,66,1259,2389\}$ or $t$ satisfies that $t \not \equiv 3,21 \bmod 27$ and that $v_p(\Delta_t) \not\equiv 2 \bmod 3$ for all $p \neq 3$.  
\end{enumerate}
In such a case, a power integral basis is given by {\rm(\ref{pib})}.
\end{crl}

The outline of this paper is as follows. In \S \ref{exm}, we present examples of \cref{thm1}, \cref{crl}.
In \S \ref{ramify}, we recall well-known properties of the conductors of cubic cyclic fields and simplest cubic fields.
The main results are in \S \ref{mvss}, \ref{cvst}.
In \S \ref{mvss}, we first show that  
\begin{align*}
&\text{$K$ has a power integral basis} &&\Longrightarrow &&\text{$\mathfrak C_K=(\alpha)$ with $\mathit{Tr}(\alpha)=0$} && \text{(\cref{fortriv})},  \\
&\text{$\mathfrak C_K=(\alpha)$ with $\mathit{Tr}(\alpha)=0$}&& \Longrightarrow &&\text{$K$ is a simplest cubic field} &&\text{(\cref{alt})}. 
\end{align*}
Hence we obtain \cref{thm2}-(i):
\begin{center}
$K$ has a power integral basis $\Lr$ 
$\begin{cases}
\text{$K$ is a simplest cubic field}, \\[3pt]
\text{$\mathfrak C_K$ is principal}. 
\end{cases}$
\end{center}
The ``converse'' does not hold true:
in the remaining of \S \ref{mvss}, we prove \cref{thm2}-(ii) which states that 
\begin{center}
$K$ has a power integral basis $\Llr$ 
$\begin{cases}
\text{$K$ is a simplest cubic field} \\[3pt]
\text{$\mathfrak C_K$ is principal} \\[3pt]
\text{one of (b), (c), (d) in \cref{thm2}-(ii)}.
\end{cases}$.
\end{center}
One of the key ideas is a relation between 
\begin{center}
``power integral bases'' \quad and \quad ``units $u$ satisfying $1+u+u\sigma(u)=0$'',
\end{center}
which can be obtained by combining \cref{fortriv}-(ii)-[(a) $\LR$ (c)] with \cref{alt}-[(i) $\LR$ (ii)].
In \S \ref{cvst}, we study the conditions when the ideal $\mathfrak C_{K_t}$ is principal and give a proof of \cref{thm3}.
Summarizing the above, we derive \cref{thm1} in \S \ref{pf}.

\begin{rmk}
\begin{enumerate}
\item Our basic idea is to decompose the monogenity of $K$ into two steps:
\begin{itemize}
\item $\mathfrak C_K$ is principal where one of its generator $\alpha$ satisfying $\mathit{Tr}(\alpha)=0$.   
\item There exists $\gamma \in \mathcal O_K$ satisfying $\alpha=\gamma-\sigma(\gamma)$ (i.e, a special case of the integral version of additive Hilbert's Theorem 90),
which have to be a power integral basis.
\end{itemize}
This idea was established in {\rm \cite{Se}} by the second author.
\item Some arguments in this paper can be seen also in {\rm\cite{Cu}}. Moreover, {\rm\cite[Lemma 1]{Cu}} states that 
if $\Delta_t$ is square-free then $K_t$ has a power integral basis.
This is a special case of {\rm\cref{thm1}-[(iii) $\R$ (i)]} in this paper.
\item Let $K$ be a cyclic cubic field. Gras {\rm\cite[Th\'eor\`eme 2]{Gr1}} showed that the followings are equivalent.
\begin{enumerate}
\item $K$ has a power integral basis.
\item There exists a unit $u \in \mathcal O_K^\times$ satisfying 
\begin{align*}
\mathit{N}(u)=1,  && \mathit{Tr}(u+u^{-1})+3=0, && \mathit{Tr} \left(\frac{u^2-u^{-1}}{c_K}\right) \in \mathbb N^3.
\end{align*}
\end{enumerate}
This corresponds to {\rm\cref{thm1}-[(i) $\LR$ (ii)]} in this paper and we provide an alternative proof.
In fact, that $\mathit{N}(u)=1$, $\mathit{Tr}(u+u^{-1})+3=0$ implies $\mathit{Tr}(u)+\mathit{Tr}(uu')+3=0$ where $u'$ is a conjugate of $u$.
It follows that $u$ is a root of $f_t(x)$ with $t:=\mathit{Tr}(u)$, i.e., $K=K_t$.
Then we see that $\mathit{Tr}(u^2-u^{-1})=\mathit{Tr}(u)^2-3\mathit{Tr}(uu')=t^2+3t+9$. 
\end{enumerate}
\end{rmk}

\section*{Acknowledgements}

The authors would like to thank Professors Toru Nakahara, Ryotaro Okazaki and Hiroshi Tsunogai for useful discussions.
We also thank Master's student Yudai Tanaka for suggesting to apply some techniques in ``genus theory'' for this topic.

\section{Examples} \label{exm}

The left table in \cref{-1to70} below is a list of $t$, the prime factorization of $c_{K_t}$ (the ramified primes), the prime factorization of $\Delta_t$, 
the corresponding ``case'', and the other $t'$ satisfying $K_t=K_{t'}$. There are 3 ``cases'':
\begin{enumerate}
\item[(a)] $K_t$ has a power integral basis by \cref{thm1} (or \cref{crl}).
\item[(b)] $K_t$ has a power integral basis since another $t'$ with $K_t=K_{t'}$ satisfies the condition of \cref{thm1}
(although $t$ does not satisfy the condition). 
\item[(c)] $K_t$ does not have a power integral basis by \cref{crl}.
\end{enumerate}
For example, 
\begin{itemize}
\item $K_{-1},K_0,K_1,K_2$ have power integral bases since these satisfy \cref{thm1}-(ii), (iii).
\item $K_{3}$ has a power integral basis since $K_{3}=K_0$.
\item $K_{21}$ is the first example which does not have a power integral basis (by \cref{crl}).
\end{itemize}
The right table is a list of $K_t$ which do \underline{\textit{not}} have power integral bases up to $t=2000$.

\begin{table}[h]
\centering 
\scriptsize
\begin{tabular}{|c|ccc|ccc|c|c|} \hline
$t $&\multicolumn{3}{|c|}{$c_{K_t}$} &\multicolumn{3}{|c|}{$\Delta_t$} & case &$=K_{t'}$\\ \hline
-1&7& &  &$7^1 $&$ $&$ $& (a)& 5,12,1259 \\ \hline
0&$3^2$&  &  &$3^2 $&$ $&$ $& (a)&3,54\\ \hline
1&13&  &  &$13^1 $&$ $&$ $& (a)&66\\ \hline
2&19&  &  &$19^1 $&$ $&$ $& (a)&2389\\ \hline
3&$3^2$&  &  &$3^3 $&$ $&$ $& (b)&0,54\\ \hline
4&37&  &  &$37^1 $&$ $&$ $& (a)&\\ \hline
5&7&  &  &$7^2 $&$ $&$ $& (b)&-1,12,1259\\ \hline
6&$3^2$&7&  &$3^2 $&$ 7^1 $&$ $& (a)&\\ \hline
7&79&  &  &$79^1 $&$ $&$ $& (a)&\\ \hline
8&97&  &  &$97^1 $&$ $&$ $& (a)&\\ \hline
9&$3^2$&13&  &$3^2 $&$ 13^1 $&$ $& (a)&\\ \hline
10&139&  &  &$139^1 $&$ $&$ $& (a)&\\ \hline
11&163&  &  &$163^1 $&$ $&$ $& (a)&\\ \hline
12&7&  &  &$3^3 $&$ 7^1 $&$ $& (a)&-1,5,1259\\ \hline
13&7&31&  &$7^1 $&$ 31^1 $&$ $& (a)&\\ \hline
14&13&19&  &$13^1 $&$ 19^1 $&$ $& (a)&\\ \hline
15&$3^2$&31&  &$3^2 $&$ 31^1 $&$ $& (a)&\\ \hline
16&313&  &  &$313^1 $&$ $&$ $& (a)&\\ \hline
17&349&  &  &$349^1 $&$ $&$ $& (a)&\\ \hline
18&$3^2$&43&  &$3^2 $&$ 43^1 $&$ $& (a)&\\ \hline
19&7&61&  &$7^1 $&$ 61^1 $&$ $& (a)&\\ \hline
20&7&67&  &$7^1 $&$ 67^1 $&$ $& (a)&\\ \hline
21&$3^2$&19&  &$3^3 $&$ 19^1 $&$ $&(c)&\\ \hline
22&13&43&  &$13^1 $&$ 43^1 $&$ $& (a)&\\ \hline
23&607&  &  &$607^1 $&$ $&$ $& (a)&\\ \hline
24&$3^2$&73&  &$3^2 $&$ 73^1 $&$ $& (a)&\\ \hline
25&709&  &  &$709^1 $&$ $&$ $& (a)&\\ \hline
26&7&109&  &$7^1 $&$ 109^1 $&$ $& (a)&\\ \hline
27&$3^2$&7&13&$3^2 $&$ 7^1 $&$ 13^1 $& (a)&\\ \hline
28&877&  &  &$877^1 $&$ $&$ $& (a)&\\ \hline
29&937&  &  &$937^1 $&$ $&$ $& (a)&\\ \hline
30&$3^2$&37&  &$3^3 $&$ 37^1 $&$ $&(c)&\\ \hline
\end{tabular} \quad
\begin{tabular}{|l|} \hline
\hspace{28pt} $K_t$ ($-1\leqq t\leqq 2000$) of case (c) \\ \hline
21,
30,
41,
48,
57,
75,
84,
90,
100,
102,
103,
111, \\ \hline
129,
138,
139,
152,
154,
156,
165,
183,
188,
192, \\ \hline
201,
204,
210,
219,
235,
237,
246,
250,
264,
269, \\ \hline 
271,
273,
291,
299,
300,
318,
327,
335,
345,
348, \\ \hline 
354,
356,
372,
374,
381,
384,
398,
399,
404,
408, \\ \hline 
426,
433,
435,
438,
446,
453,
462,
480,
482,
489, \\ \hline 
495,
507,
515,
516,
531,
534,
543,
544,
561,
565, \\ \hline 
570,
573,
577,
580,
588,
593,
597,
602,
607,
615, \\ \hline 
624,
642,
651,
669,
678,
691,
696,
705,
716,
723, \\ \hline 
727,
732,
742,
750,
759,
776,
777,
786,
789,
804, \\ \hline 
813,
825,
831,
838,
840,
844,
858,
867,
874,
876, \\ \hline 
885,
887,
894,
912,
921,
923,
926,
936,
939,
945, \\ \hline 
948,
966,
975,
985,
992,
993,
1002,
1020,
1021, \\ \hline 
1029,
1034,
1047,
1056,
1070,
1074,
1080,
1083, \\ \hline 
1096,
1101,
1110,
1114,
1119,
1128,
1132,
1137, \\ \hline 
1155,
1164,
1168,
1181,
1182,
1191,
1209,
1210, \\ \hline 
1217,
1218,
1230,
1236,
1237,
1245,
1249,
1263, \\ \hline 
1265,
1266,
1269,
1272,
1279,
1287,
1290,
1299, \\ \hline 
1317,
1326,
1328,
1344,
1353,
1364,
1371,
1377, \\ \hline 
1380,
1398,
1407,
1413,
1418,
1425,
1434,
1452, \\ \hline 
1461,
1462,
1475,
1479,
1488,
1497,
1506,
1511, \\ \hline 
1515,
1524,
1533,
1542,
1560,
1563,
1569,
1573, \\ \hline 
1587,
1596,
1609,
1614,
1621,
1622,
1623,
1641, \\ \hline 
1648,
1650,
1668,
1671,
1677,
1695,
1704,
1707, \\ \hline 
1720,
1722,
1731,
1743,
1749,
1756,
1758,
1776, \\ \hline 
1785,
1790,
1803,
1805,
1812,
1818,
1830,
1839, \\ \hline 
1854,
1857,
1866,
1867,
1884,
1893,
1903,
1911, \\ \hline 
1916,
1920,
1925,
1938,
1947,
1952,
1959,
1965, \\ \hline 
1974,
1992 \\ \hline 
\end{tabular}
\caption{small $t$}  \label{-1to70}
\end{table}

\newpage

\cref{nontriv} below is a list of ``non-trivial cases'': those satisfy $\mathcal O_{K_t} \supsetneqq \mathbb Z[\theta_t]$ and $\mathfrak C_{K_t}$ is principal.
Note that 
\begin{itemize}
\item When $\mathcal O_{K_t} = \mathbb Z[\theta_t]$ (i.e., $c_{K_t}=\Delta_t$), it has a (trivial) power integral basis $\theta_t$.
\item When $\mathfrak C_{K_t}$ is not principal, it does not have a power integral basis by \cref{thm2}.
\item If $t< 101471$ then the converse of {\rm\cref{thm2}-(i)}: ``$\mathfrak C_{K_t}$ is principal $\R$ $K_t$ has a power integral basis'' also holds true.
All the examples of 
\begin{center}
``$\mathfrak C_{K_t}$ is principal but $K_t$ does not have a power integral basis'', \\
up to $t=10^8=100000000$
\end{center}
are $t=101471,182451,18128865$. 
\item $K_{740}$ has a power integral basis although $v_7(\Delta_{740})=4 \neq 1$. 
This is the first example of ``$K_t$ has a power integral basis although the prime-to-$3$ part of $\Delta_t$ is not square-free'',
except for the exceptions $K_5$ (=$K_{-1}$), $K_{54}$ (=$K_0$), $K_{66}$ (=$K_1$). 
\end{itemize}

\begin{table}[h]
\centering 
\scriptsize
\begin{tabular}{|c|ccc|ccccc|c|c|c|} \hline
$t $&\multicolumn{3}{|c|}{$c_{K_t}$}&\multicolumn{5}{|c|}{$\Delta_t$} & case &$=K_{t'}$\\ \hline
3&$3^2$ &  &  &$3^3 $&$ $&$ $&$ $&$ $& (b)&0,54\\ \hline
5&7 &  &  &$7^2 $&$ $&$ $&$ $&$ $& (b)&-1,12,1259\\ \hline
12&7 &  &  &$3^3 $&$ 7^1 $&$ $&$ $&$ $& (a)&-1,5,1259\\ \hline
39&61 &  &  &$3^3 $&$ 61^1 $&$ $&$ $&$ $& (a)&\\ \hline
54&$3^2$ &  &  &$3^2 $&$ 7^3 $&$ $&$ $&$ $& (a)&0,3\\ \hline
66&13 &  &  &$3^3 $&$ 13^2 $&$ $&$ $&$ $& (b)&1\\ \hline
93&331 &  &  &$3^3 $&$ 331^1 $&$ $&$ $&$ $& (a)&\\ \hline
120&547 &  &  &$3^3 $&$ 547^1 $&$ $&$ $&$ $& (a)&\\ \hline
147&19 & 43 &  &$3^3 $&$ 19^1 $&$ 43^1 $&$ $&$ $& (a)&\\ \hline
174&7 & 163 &  &$3^3 $&$ 7^1 $&$ 163^1 $&$ $&$ $& (a)&\\ \hline
228&1951 &  &  &$3^3 $&$ 1951^1 $&$ $&$ $&$ $& (a)&\\ \hline
255&2437 &  &  &$3^3 $&$ 2437^1 $&$ $&$ $&$ $& (a)&\\ \hline
282&13 & 229 &  &$3^3 $&$ 13^1 $&$ 229^1 $&$ $&$ $& (a)&\\ \hline
286&241 &  &  &$7^3 $&$ 241^1 $&$ $&$ $&$ $& (a)&\\ \hline
309&3571 &  &  &$3^3 $&$ 3571^1 $&$ $&$ $&$ $& (a)&\\ \hline
336&4219 &  &  &$3^3 $&$ 4219^1 $&$ $&$ $&$ $& (a)&\\ \hline
363&7 & 19 & 37 &$3^3 $&$ 7^1 $&$ 19^1 $&$ 37^1 $&$ $& (a)&\\ \hline
390&7 & 811 &  &$3^3 $&$ 7^1 $&$ 811^1 $&$ $&$ $& (a)&\\ \hline
397&463 &  &  &$7^3 $&$ 463^1 $&$ $&$ $&$ $& (a)&\\ \hline
417&13 & 499 &  &$3^3 $&$ 13^1 $&$ 499^1 $&$ $&$ $& (a)&\\ \hline
444&7351 &  &  &$3^3 $&$ 7351^1 $&$ $&$ $&$ $& (a)&\\ \hline
471&8269 &  &  &$3^3 $&$ 8269^1 $&$ $&$ $&$ $& (a)&\\ \hline
498&9241 &  &  &$3^3 $&$ 9241^1 $&$ $&$ $&$ $&(a)&\\ \hline
525&10267 &  &  &$3^3 $&$ 10267^1 $&$ $&$ $&$ $&(a)&\\ \hline
552&7 & 1621 &  &$3^3 $&$ 7^1 $&$ 1621^1 $&$ $&$ $&(a)&\\ \hline
579&7 & 1783 &  &$3^3 $&$ 7^1 $&$ 1783^1 $&$ $&$ $&(a)&\\ \hline
606&13669 &  &  &$3^3 $&$ 13669^1 $&$ $&$ $&$ $&(a)&\\ \hline
629&19 & 61 &  &$7^3 $&$ 19^1 $&$ 61^1 $&$ $&$ $&(a)&\\ \hline
633&13 & 31 & 37 &$3^3 $&$ 13^1 $&$ 31^1 $&$ 37^1 $&$ $&(a)&\\ \hline
660&19 & 853 &  &$3^3 $&$ 19^1 $&$ 853^1 $&$ $&$ $&(a)&\\ \hline
687&97 & 181 &  &$3^3 $&$ 97^1 $&$ 181^1 $&$ $&$ $&(a)&\\ \hline
714&67 & 283 &  &$3^3 $&$ 67^1 $&$ 283^1 $&$ $&$ $&(a)&\\ \hline
740&7 & 229 &  &$7^4 $&$ 229^1 $&$ $&$ $&$ $& (a)&\\ \hline 
\hline
101471&5479 &  &  &$7^3 $&$ 5479^2 $&$ $&$ $&$ $& (c)&  \\ \hline \hline
182451&13 & 37 & 73 &$3^3 $&$ 13^2 $&$ 37^2 $&$ 73^2 $&$ $& (c)&  \\ \hline \hline
18128865&13 & 43 & 337 &$3^3 $&$ 7^3 $&$ 13^2 $&$ 43^2 $&$ 337^2 $& (c)&  \\ \hline
\end{tabular}
\caption{non-trivial ($\mathcal O_{K_t} \supsetneq \mathbb Z[\theta_t]$, $\mathfrak C_{K_t}$: principal) cases} \label{nontriv}
\end{table}

\section{Primes dividing $c_K,\Delta_t$} \label{ramify}

In this section, we recall properties of primes dividing $c_K$ or $\Delta_t$ which are well-known for experts.
First we note that 
\begin{align*}
\{\text{all the ramified primes}\} = \{p \mid c_{K_t}\} \subset \{p \mid \Delta_t\}
\end{align*}
since $c_{K_t} \mid \Delta_t$. 
The following Propositions are well-known for experts.

\begin{prp} \label{order}
Let $K$ be a cyclic cubic field. The conductor $c_K$ satisfies the following conditions. Let $p$ denote a rational prime.
\begin{enumerate}
\item If $3 \mid c_K$, then $v_3(c_K)=2$.
\item If $p \equiv 1 \bmod 3$, then  $p \mid c_K$ implies $v_p(c_K)=1$.
\item If $p \equiv 2 \bmod 3$, then $p \nmid c_K$.
\end{enumerate}
\end{prp}

\begin{prp} \label{3}
Let $K=K_t$. Then $v_3(\Delta_t)$ takes only values of $0,2,3$. More precisely, we have the following. 
\begin{enumerate}
\item $v_3(\Delta_t)=0$ if and only if $3 \nmid t$. In this case, $3 \nmid c_K$.
\item $v_3(\Delta_t)=2$ if and only if $t \equiv 0,6 \bmod 9$. In this case, $3 \mid c_K$.
\item $v_3(\Delta_t)=3$ if and only if $t \equiv 3 \bmod 9$. In this case, 
\begin{align*}
\begin{cases}
t\equiv 12 \bmod 27 &\Llr 3 \nmid c_K, \\
t\equiv 3,21 \bmod 27 &\Llr 3 \mid c_K.
\end{cases}
\end{align*}
In particular we see that 
\begin{align} \label{when}
\begin{cases}
3 \nmid c_K \Llr 3 \nmid t \text{ or } t\equiv 12 \bmod 27, \\
3 \mid c_K \Llr t \equiv 0,6 \bmod 9 \text{ or } t\equiv 3,21 \bmod 27.
\end{cases}
\end{align}
\end{enumerate}
\end{prp}

\begin{prp} \label{not3}
Let $K=K_t$, $p \neq 3$. 
The followings are equivalent.
\begin{align*}
v_p(\Delta_t) \equiv 0 \bmod 3 \Llr p \nmid c_K.
\end{align*}
\end{prp}

\cref{order} can be shown by noting that $c_K$ is the minimal integer satisfying $K_t \subset \mathbb Q(\zeta_{c_K})$, which implies 
\begin{align*}
\prod_{p \mid c_K} (\mathbb Z/p^{v_p(c_K)}\mathbb Z)^\times\cong (\mathbb Z/c_K\mathbb Z)^\times \cong \mathrm{Gal}(\mathbb Q(\zeta_{c_K})/\mathbb Q) 
\twoheadrightarrow \mathrm{Gal}(\mathbb Q(K/\mathbb Q) \cong \mathbb Z/3\mathbb Z.
\end{align*}
Here, we give a proof only of \cref{3} since that of \cref{not3} is similar and simpler.

\begin{proof}[Proof of \rm{\cref{3}}]
The ``if and only if'' part follows from an explicit calculation as  
\begin{align*}
&\Delta_{3t_0}= 9t_0^2 + 9t_0 + 9, && \Delta_{3t_0+1}= 9t_0^2 + 15t_0 + 13, && \Delta_{3t_0+2}=9t_0^2 + 21t_0 + 19, \\
&\Delta_{9t_0}=9(9t_0^2 + 3t_0 + 1), && \Delta_{9t_0+3}=27( 3t_0^2 + 3t_0 + 1), && \Delta_{9t_0+6}=9(9t_0^2 + 15t_0 + 7).
\end{align*}
Then (i)  holds true since $c_K \mid \Delta_t$. \\[5pt]
(ii) Let $t=9t_0$. Then the minimal polynomial of $\theta_t-\frac{t}{3}-1$:
\begin{align*}
x^3 + 3x^2 + (-27t_0^2 - 9t_0)x + (-54t_0^3 - 54t_0^2 - 18t_0 - 3)
\end{align*}
is an Eisenstein polynomial, so we have $3\mid c_K$. Similarly, let $t=9t_0+6$. Then the minimal polynomial of $\theta_t-\frac{t}{3}+1$:
\begin{align*}
x^3 - 3x^2 + (-27t_0^2 - 45t_0 - 18)x + (-54t_0^3 - 108t_0^2 - 72t_0 - 15)
\end{align*}
also is an Eisenstein polynomial, so we have $3\mid c_K$. \\[5pt] 
(iii) Assume $t\equiv 3 \bmod 9$. Then the minimal polynomial of $\theta_t-\frac{t}{3}$ is
\begin{align*} 
x^3-\frac{\Delta_t}{3}x-\frac{(2t+3)\Delta_t}{3^3} \in \mathbb Z[x].
\end{align*}
Here we note that $v_3(\Delta_t)=3$.
First assume that $t\equiv 3,21 \bmod 27$, which implies $v_3(2t+3 )=2$. Then the Newton polygon is given by
\begin{align*}
v_3(1)=0, && v_3(0)=\infty, && v_3(\tfrac{\Delta_t}{3})=2, && v_3(\tfrac{(2t+3)\Delta_t}{3^3})= 2,
\end{align*}
so we have $v_3(\theta_t-\frac{t}{3})=\frac{2}{3}$. That is, $3$ ramifies in $K_t=\mathbb Q(\theta_t-\frac{t}{3})$.
Next assume that $t\equiv 12 \bmod 27$, which implies $v_3(2t+3 )\geqq 3$.
In this case, the Newton polygon is given by
\begin{align*}
v_3(1)=0, && v_3(0)=\infty, && v_3(\tfrac{\Delta_t}{3})=2, && v_3(\tfrac{(2t+3)\Delta_t}{27})\geqq 3.
\end{align*}
That is, $3 \mid \theta_t-\frac{t}{3}$. Then the discriminant of (the minimal polynomial of) $\frac{\theta_t-\frac{t}{3}}{3}$
equals $(\frac{\Delta_t}{3^3})^2$, which is prime to $3$. 
Note that $c_K \mid \frac{\Delta_t}{3^3}$ since $\mathbb Z\left[\frac{\theta_t-\frac{t}{3}}{3}\right] \subset \mathcal O_{K_t}$.
Then we obtain $3 \nmid c_K$ as desired.
\end{proof}

\begin{proof}[Proof of the expression {\rm(\ref{factor})}]
The division into two cases comes from whether $3 \mid c_{K_t}$ or not, by \cref{3}.
The range of $p$ and its exponent follow from \cref{not3} and \cref{order}, respectively.
\end{proof}

\section{Monogenity and simplest cubic fields} \label{mvss}

In this section, we give a proof of \cref{thm2} which provides an equivalent condition for that a cyclic cubic field $K$ has a power integral basis. 

\begin{prp} \label{fortriv}
Let $K$ be a cyclic cubic field.
\begin{enumerate}
\item There exists a unique integral ideal $\mathfrak C_K \subset \mathcal O_K$ satisfying $\mathit{N} \mathfrak C_K=c_K$.
\item The followings are equivalent for $\gamma \in \mathcal O_K$.
\begin{enumerate}
\item $\gamma$ is a power integral basis. That is, $\mathcal O_K=\mathbb Z[\gamma]$.
\item $c_K=\mathit{N}(\gamma-\sigma(\gamma))$.
\item $\mathfrak C_K=(\gamma-\sigma(\gamma))$.
\end{enumerate}
\end{enumerate}
In particular, if $K$ has a power integral basis, then $\mathfrak C_K$ is a principal ideal with a generator $\alpha:=\gamma-\sigma(\gamma)$ satisfying $\mathit{Tr}(\alpha)=0$. 
\end{prp}

\begin{proof}
(i) The followings are equivalent
\begin{center}
$p \mid c_K \Llr p\mid d_K \Llr$ a unique prime ideal $\mathfrak P_p$ satisfies $p=\mathfrak P_p^3$, $\mathit{N}\mathfrak P_p=p$. 
\end{center}
Hence 
\begin{align*}
\mathfrak C_K:=\prod_{p \mid c_K} \mathfrak P_p^{v_p(c_K)}
\end{align*}
is the unique ideal satisfying $\mathit{N} \mathfrak C_K=c_K$. \\[5pt]
(ii) We put $d(\gamma):=\prod_{\iota \neq \iota'}(\iota(\gamma)-\iota'(\gamma))=\mathit{N}(\gamma-\sigma(\gamma))^2$,
where $\iota,\iota'$ run over all embeddings of $K$ with $\iota \neq \iota'$. Then we can write
\begin{align*}
\mathcal O_K=\mathbb Z[\gamma] \Llr d_K=d(\gamma) \Llr c_K=\mathit{N}(\gamma-\sigma(\gamma)),
\end{align*}
so the assertion follows from (i).
\end{proof}

\begin{lmm} \label{alt}
Let $K$ be a cyclic cubic field. Assume that $\mathfrak C_K$ is principal, take a generator $\beta$ of $\mathfrak C_K$, 
and put $u_\beta:=\beta^{\sigma-1} \in \mathcal O_K^\times$.
The followings are equivalent.
\begin{enumerate}
\item There exists $\alpha \in \mathcal O_K$ satisfying $\mathfrak C_K=(\alpha)$, $\mathit{Tr}(\alpha)=0$.
\item There exists $u \in \mathcal O_K^\times$ satisfying $1+u+u\sigma(u)=0$, $[u_\beta]=[u] \in H^1(G,\mathcal O_K^\times)$.
\item There exists $t$ satisfying $K=K_t$, $[u_\beta]=[\theta_t] \in H^1(G,\mathcal O_K^\times)$.
\end{enumerate}
\end{lmm}

\begin{proof}
{[(i) $\LR$ (ii)]}. For any generator $\alpha:=\beta\epsilon \in \mathfrak C_K$ $(\epsilon \in \mathcal O_K^\times)$, we have
\begin{align*}
\mathit{Tr}(\alpha)=\beta\epsilon+\sigma(\beta\epsilon)+\sigma^2(\beta\epsilon)=\beta\epsilon(1+u_\beta\epsilon^{\sigma-1} +u_\beta\epsilon^{\sigma-1}\sigma(u_\beta\epsilon^{\sigma-1})).
\end{align*}
In particular, the followings are equivalent as desired: 
\begin{center}
$\mathit{Tr}(\alpha)=0$ $\Llr$ $\exists u:=u_\beta\epsilon^{\sigma-1} \in [u_\beta]$ s.t.\ $1+u+u\sigma(u)=0$.
\end{center}
[(ii) $\LR$ (iii)]. ($\L$) is obvious. For ($\R$), we note that 
\begin{center}
if $1+u+u\sigma(u)=0$, then we have $0=\mathit{Tr}(1+u+u\sigma(u))=3+\mathit{Tr}(u)+\mathit{Tr}(u\sigma(u))$.
\end{center}
On the other hand, by $u \in \mathcal O_K$, we have
\begin{align*}
t:=\mathit{Tr}(u) \in \mathbb Z.
\end{align*}
Therefore, $u$ is a root of 
\begin{align*}
(x-u)(x-\sigma(u))(x-\sigma^2(u))&=x^3-\mathit{Tr}(u)x^2+\mathit{Tr}(u\sigma(u))x-1 =f_t(x).
\end{align*}
Hence $u$ is a conjugate $\tau(\theta_t)$ of $\theta_t$ ($\tau \in G$), so we have $u=\tau(\theta_t)=\theta_t(\theta_t)^{\tau-1} \in [\theta_t]$.
\end{proof}

\begin{proof}[Proof of {\rm\cref{thm2}-(i)}]
Let $K$ be a cyclic cubic field with a power integral basis $\gamma$. Then $\mathfrak C_K$ is principal by \cref{fortriv}.
Moreover its generator $\alpha:=\gamma-\sigma(\gamma)$ satisfies $\mathit{Tr}(\alpha)=0$.
It follows that $K$ is a simplest cubic field by \cref{alt}. 
\end{proof}

We show the remaining statements of \cref{thm2}.

\begin{prp} \label{coh}
\begin{enumerate}
\item We have a canonical isomorphism 
\begin{align*}
P_K^G/P_{\mathbb Q} \cong H^1(G,\mathcal O_K^\times), \ \overline{(\alpha)} \mt [\alpha^{\sigma-1}].
\end{align*}
\item The norm map induces the following injective homomorphism 
\begin{align*}
I_K^G/P_{\mathbb Q} \hookrightarrow \mathbb Q^\times/(\mathbb Q^\times)^3, \ \overline{\mathfrak A} 
\mt \mathit{N}\mathfrak A \bmod (\mathbb Q^\times)^3.
\end{align*}
\end{enumerate}
In particular, we have
\begin{align*}
H^1(G,\mathcal O_K^\times) \cong P_K^G/P_{\mathbb Q} \subset I_K^G/P_{\mathbb Q} \hookrightarrow \mathbb Q^\times/(\mathbb Q^\times)^3.
\end{align*}
\end{prp}

\begin{proof}
(i) Consider the exact sequence 
\begin{align*}
1 \ra \mathcal O_K^\times \ra K^\times \ra P_K \ra 1.
\end{align*}
Then we have by Hilbert's Theorem 90
\begin{align*}
1 \ra \{\pm 1\} \ra \mathbb Q^\times \ra P_K^G \ra H^1(G,\mathcal O_K^\times) \ra 1. 
\end{align*}
Hence we get 
\begin{align*}
H^1(G,\mathcal O_K^\times) \cong [\mathrm{Coker}\colon \mathbb Q^\times \ra P_K^G] =P_K^G/P_{\mathbb Q},\  [\alpha^{\sigma-1}] \mapsfrom \overline{(\alpha)}
\end{align*}
by calculating the connecting homomorphism. \\[5pt]
(ii) Since $I_K^G$ is generated by ramified prime ideals and ideals contained in $P_\mathbb Q$, we can write any ideal contained in $ I_K^G$ as
\begin{align*}
\mathfrak A:=(n)\prod_{p \mid c_K} \mathfrak P_p^{e_p} \text{ with }\mathfrak P_p^3=(p),\ n \in \mathbb N. 
\end{align*}
Then the followings are equivalent:
\begin{align*}
\mathit{N} \mathfrak A=(n^3)\prod_{p \mid c_K} p^{e_p} \in (\mathbb Q^\times)^3 \Llr e_p \in 3\mathbb Z \Llr \mathfrak A \in P_\mathbb Q.
\end{align*}
Hence the kernel of the norm map $I_K^G \ra \mathbb Q^\times/(\mathbb Q^\times)^3$ equals $P_\mathbb Q$.
Then the assertion is clear.
\end{proof}

\begin{proof}[Proof of {\rm\cref{thm2}-(ii)}]
Assume that $K$ is a simplest cubic field and that $\mathfrak C_K=(\beta)$. 
The equivalences (b) $\LR$ (c) $\LR$ (d) follow from \cref{coh} by noting that 
\begin{align*}
\begin{array}{ccccc}
H^1(G,\mathcal O_K^\times) &\cong& P_K^G/P_{\mathbb Q} &\hookrightarrow& \mathbb Q^\times/(\mathbb Q^\times)^3, \\
\text{\rotatebox[origin=c]{90}{$\in$}} && \text{\rotatebox[origin=c]{90}{$\in$}} && \text{\rotatebox[origin=c]{90}{$\in$}} \\
{[u_\beta]} & \mapsfrom & \overline{\mathfrak C_K} & \mt & c_K \bmod (\mathbb Q^\times)^3, \\[3pt]
{[\theta_t]} & \mapsfrom & \overline{(\theta_t-\sigma(\theta_t))} & \mt & \Delta_t \bmod (\mathbb Q^\times)^3.
\end{array}
\end{align*}
Here, by (\ref{amb}), we obtain the bottom-left part 
as $[(\theta_t-\sigma(\theta_t))^{\sigma-1}]=[\sigma^2(\theta_t)]=[\theta_t\theta_t^{\sigma^2-1}]=[\theta_t] \in H^1(G,\mathcal O_K^\times)$. 
The part ``(a) $\R$ (b)'' follows by using \cref{fortriv}-(ii)-[(a) $\R$ (c)] and \cref{alt}-[(i) $\R$ (iii)] for $\alpha:=\gamma-\sigma(\gamma)$.
We prove the remaining part ``(d) $\R$ (a)'' by showing that $\gamma$ in (\ref{pib}) actually provides a power integral basis. 
It suffices to show that
\begin{align*}
\gamma = \frac{\theta_t - a}{\sqrt[3]{\frac{\Delta_t}{c_K}}} \in \mathcal{O}_K
\end{align*}
since we have clearly $d(\gamma)=\frac{d(\theta_t)}{\left(\frac{\Delta_t}{c_K}\right)^2}=c_K^2=d_{K}$ for $a \in \mathbb Z$. \\[5pt]
[(d) $\R$ (a)]. Since $(\theta_t-\sigma(\theta_t))$ is an ambiguous ideal, we see that $(\theta_t-\sigma(\theta_t))^3=(\mathit{N}(\theta_t-\sigma(\theta_t)))$.
On the other hand, we have $\mathit{N}(\theta_t-\sigma(\theta_t)) = \Delta_t \in c_K \cdot \mathbb{N}^3$. It follows that 
\begin{align}
\theta_t \equiv \sigma(\theta_t) \mod \sqrt[3]{\frac{\Delta_t}{c_K}}. 
\label{Mod}
\end{align}
Assume $3 \nmid \sqrt[3]{\frac{\Delta_t}{c_K}}$. Then there exists $a \in \mathbb{Z}$ with $a \equiv \frac{t}{3} \bmod \sqrt[3]{\frac{\Delta_t}{c_K}}$ 
by $3 \in \left(\mathbb{Z} / \sqrt[3]{\frac{\Delta_t}{c_K}} \mathbb{Z}\right)^\times$, so we have
\begin{align*}
a \equiv \frac{t}{3} = \frac{\mathit{Tr}(\theta_t) }{3} \equiv \frac{3\theta_t}{3} = \theta_t \mod \sqrt[3]{\frac{\Delta_t}{c_K}}
\end{align*}
by (\ref{Mod}). Hence we have $\theta_t - a \in \sqrt[3]{\frac{\Delta_t}{c_K}} \mathcal{O}_K$, namely $\gamma = \frac{\theta_t - a}{\sqrt[3]{\frac{\Delta_t}{c_K}}} \in \mathcal{O}_K$. 

Assume $3 \mid \sqrt[3]{\frac{\Delta_t}{c_K}}$. Then we obtain $v_3(\Delta_t)=3$, $3 \nmid c_K$, $v_3 \left( \sqrt[3]{\frac{\Delta_t}{c_K}} \right) = 1$, 
$t\equiv 12 \bmod 27$ by \cref{3}.
In this case, by (\ref{Mod}), we have
\begin{align*}
t = \mathit{Tr}(\theta_t) \equiv 3\theta_t  \mod \sqrt[3]{\frac{\Delta_t}{c_K}}.
\end{align*}
It follows that  
\begin{align} \label{add}
\theta_t - \frac{t}{3} \in \frac{1}{3} \sqrt[3]{\frac{\Delta_t}{c_K}} \mathcal{O}_K.
\end{align}
On the other hand, in Proof of \cref{3}-(iii) in \S \ref{ramify}, we showed that $3 \mid \theta_t - \frac{t}{3}$ when $t\equiv 12 \bmod 27$.
By combining this with (\ref{add}), we obtain $\gamma = \frac{\theta_t - a}{\sqrt[3]{\frac{\Delta_t}{c_K}}} \in \mathcal{O}_K$ for $a=\frac{t}{3}$.
\end{proof}

\section{Relation between $\mathfrak C_{K_t}$ and $(\theta_t-\sigma(\theta_t))$} \label{cvst}

Let $K=K_t$ be a simplest cubic field generated by a root $\theta_t$ of $f_t(x)$. Recall that 
\begin{itemize}
\item The ambiguous ideal $\mathfrak C_{K}$ satisfies $\mathit{N} \mathfrak C_{K}=c_{K}$. 
\item The principal ambiguous ideal $(\theta_t-\sigma(\theta_t))$ satisfies $\mathit{N}(\theta_t-\sigma(\theta_t))=\Delta_t$.
\item By \cref{coh} and Proof of \cref{thm2}-(ii) in the previous section, we have  
\begin{align} \label{diag}
\begin{array}{ccccccc}
H^1(G,\mathcal O_K^\times) &\cong& P_K^G/P_{\mathbb Q} &\subset& I_K^G/P_{\mathbb Q}  &\hookrightarrow& \mathbb Q^\times/(\mathbb Q^\times)^3, \\
\text{\rotatebox[origin=c]{90}{$\in$}} && \text{\rotatebox[origin=c]{90}{$\in$}} && \text{\rotatebox[origin=c]{90}{$\in$}} \\
& & & & \overline{\mathfrak C_K} & \mt & c_K \bmod (\mathbb Q^\times)^3, \\[3pt]
{[\theta_t]} & \mapsfrom & \overline{(\theta_t-\sigma(\theta_t))} & & & \mt & \Delta_t \bmod (\mathbb Q^\times)^3
\end{array}
\end{align}
without assuming that $\mathfrak C_K$ is principal.
\end{itemize}

\begin{prp} \label{mu}
Let $K$ be a cyclic cubic field. We have
\begin{align*}
H^1(G,\mathcal O_K^\times) \cong \mathbb Z/3\mathbb Z.
\end{align*}
\end{prp}

\begin{proof}
For a finite cyclic extension $L/F$, the Herbrand quotient $Q(\mathcal O_L^\times)$ equals 
$\frac{1}{[L:F]} \times $ ``the product of ramification indices of all the infinite places'' (for a proof, see \cite[Lemma 3]{Yo}).
Therefore we see that 
\begin{align*}
\frac{1}{|H^1(G,\mathcal O_K^\times) |}=Q(\mathcal O_K^\times)=\frac{1}{3},
\end{align*}
since $\mathit{N}\colon \mathcal O_K^\times \ra \mathbb Z^\times$ is surjective.
Then the assertion is clear.
\end{proof}

\begin{proof}[Proof of {\rm\cref{thm3}}]
Let $K=K_t \neq \mathbb Q(\zeta_9+\zeta_9^{-1})$.
First we show that $(\theta_t-\sigma(\theta_t)) \notin P_\mathbb Q$ (strictly speaking, $(\theta_t-\sigma(\theta_t))\notin [\mathrm{Im}\colon P_\mathbb Q \hookrightarrow P_K^G]$).
If $(\theta_t-\sigma(\theta_t)) \in P_\mathbb Q$, then $\mathit{N} (\theta_t-\sigma(\theta_t))=\Delta_t \in \mathbb N^3$, 
which implies that $p \nmid c_K$ for any $p\neq 3$ by \cref{not3}. 
It contradicts with $K\neq \mathbb Q(\zeta_9+\zeta_9^{-1})$. 
We see also that $\mathfrak C_K \notin P_\mathbb Q$ by \cref{order}.
Then we can write by (\ref{diag}) and \cref{mu}
\begin{align*}
\begin{array}{ccccc}
P_K^G/P_{\mathbb Q} &\subset& I_K^G/P_{\mathbb Q} &\hookrightarrow& \mathbb Q^\times/(\mathbb Q^\times)^3, \\
\text{\rotatebox[origin=c]{90}{$=$}} && \text{\rotatebox[origin=c]{90}{$\in$}} \\
\left\{\overline{(\theta_t-\sigma(\theta_t))}^e \mid e=0,1,2\right\}  & & \overline{\mathfrak C_K} \neq \overline{(1)}.
\end{array}
\end{align*}
Therefore the followings are equivalent:
\begin{align*}
\text{$\mathfrak C_K$ is principal} &\Llr \overline{\mathfrak C_K} \in P_K^G/P_{\mathbb Q}\\
&\Llr \text{there exists $e\in \{1,2\}$ satisfying }\mathfrak C_K \equiv (\theta_t-\sigma(\theta_t))^e \bmod P_\mathbb Q\\
&\Llr \text{there exists $e \in \{1,2\}$ satisfying }c_K \equiv \Delta_t^e \bmod (\mathbb Q^\times)^3.
\end{align*}
Then the assertion is clear.
\end{proof}

We also provide an explicit version of \cref{thm3}: the condition in \cref{expl}-(iii) below is written explicitly in terms of $t$.
\begin{crl} \label{expl}
Let $K=K_t\neq \mathbb Q(\zeta_9+\zeta_9^{-1})$.
The followings are equivalent.
\begin{enumerate}
\item $\mathfrak C_K$ is principal.
\item 
$
\begin{cases}
v_p(\Delta_t) \equiv 1 \bmod 3 \text{ for all } p \mid c_K\text{ or }v_p(\Delta_t) \equiv 2 \bmod 3 \text{ for all }p \mid c_K & (3 \nmid c_K), \\[3pt]
v_3(\Delta_t)=2,\ v_p(\Delta_t)\equiv 1 \bmod 3 \text{ for all } 3\neq p \mid c_K & (3 \mid c_K).
\end{cases}
$
\item 
$
\begin{cases}
v_p(\Delta_t) \not \equiv 1 \bmod 3 \text{ for all } p \text{ or }v_p(\Delta_t) \not \equiv 2 \bmod 3 \text{ for all }p \hspace{-15pt}\\
& (3 \nmid t \text{ or } t\equiv 12 \bmod 27), \\[5pt]
v_p(\Delta_t)\not \equiv 2 \bmod 3 \text{ for all } p\neq 3 & (t \equiv 0,6 \bmod 9).
\end{cases}
$
\end{enumerate}
\end{crl}

\begin{proof}
The equivalence (ii) $\LR$ (iii) follows from Propositions \ref{3}, \ref{not3} immediately.
We prove ``(i) $\LR$ (ii)''.
By Propositions \ref{order}, \ref{3}, \ref{not3}, we see that  
\begin{itemize}
\item When $p\neq 3$ we have
\begin{align*}
v_p(c_{K_t})=0,1, && v_p(c_{K_t})=0 \LR v_p(\Delta_t) \equiv 0 \bmod 3.
\end{align*}
\item For $p=3$ we have
\begin{align*}
v_3(\Delta_t)=0,2,3, && 3\nmid c_{K_t} \R v_3(\Delta_t)=0,3, && 3\mid c_{K_t} \LR 3^2\| c_{K_t} \R v_3(\Delta_t)=2,3.
\end{align*}
\end{itemize}
\cref{thm3} implies that 
\begin{align} \label{trans}
\text{(i) $\mathfrak C_K$ is principal} \Llr \text{$e=1$ or $2$ satisfies $v_p(\Delta_t) \equiv ev_p(c_K) \bmod 3$ for all $p$}.
\end{align}
When $3 \nmid c_K$, we may rewrite (\ref{trans}) as
\begin{align*}
&\Llr \text{$e=1$ or $2$ satisfies $v_p(\Delta_t) \equiv ev_p(c_K) \bmod 3$ for all $p\mid c_K$} \\
&\Llr \text{$e=1$ or $2$ satisfies $v_p(\Delta_t) \equiv e \bmod 3$ for all $p\mid c_K$ (that is, (ii) with $3 \nmid c_K$)},
\end{align*}
by using the facts that  
\begin{align*}
&v_p(c_K)=0 \equiv ev_p(\Delta_t) \mod 3 && (p\nmid c_K), \\
&v_p(c_K)=1 && (p\mid c_K)
\end{align*}
respectively.

When $3 \mid c_K$, we have $v_3(\Delta_t) \not\equiv 2v_p(c_K) \bmod 3$ by $v_3(\Delta_t)=2,3$, $v_3(c_K)=2$. 
Hence the case of $e=2$ in (\ref{trans}) can not happen.
Therefore we may rewrite (\ref{trans}) as
\begin{align*}
&\Llr \text{$v_p(\Delta_t) \equiv v_p(c_K) \bmod 3$ for all $p\mid c_K$} \\
&\Llr 
\begin{cases}
v_3(\Delta_t)=2, \\
v_p(\Delta_t)\equiv 1 \bmod 3 \text{ for all } 3\neq p \mid c_K
\end{cases}
\end{align*}
as desired, by using facts that $v_3(\Delta_t)=2,3$, $v_3(c_K)=2$, and $v_p(c_K)=1$ ($3\neq p\mid c_K$).
\end{proof}

\section{Proof of \cref{thm1}} \label{pf}

We give a proof of \cref{thm1}.
The expression (\ref{pib}) of a power integral basis is given in the proof of \cref{thm2}-[(d) $\R$ (a)] in \S \ref{mvss}.
We show that the following conditions in \cref{thm1} are equivalent (we rewrite \textcircled{\scriptsize iii} slightly by using \cref{3}):
\begin{enumerate}
\item[\textcircled{\scriptsize i}] $K$ has a power integral basis.
\item[\textcircled{\scriptsize ii}] There exists $t$ satisfying that $K=K_t$ and $\frac{\Delta_t}{c_K} \in \mathbb N^3$.
\item[\textcircled{\scriptsize iii}] There exists $t$ satisfying that $K=K_t$ and
\begin{align*}
\begin{cases}
v_p(\Delta_t) \not \equiv 2 \bmod 3 \text{ for all } p
& (3 \nmid t \text{ or } t\equiv 12 \bmod 27), \\[3pt]
v_p(\Delta_t) \not\equiv 2 \bmod 3 \text{ for all } p \neq 3 
& (t \equiv 0,6 \bmod 9).
\end{cases}
\end{align*}
\end{enumerate}
The equivalence \textcircled{\scriptsize i} $\LR$ \textcircled{\scriptsize ii} follows from \cref{thm2}-(i), \cref{thm2}-(ii)-[(a) $\LR$ (d)] and \cref{thm3}-[(i) $\LR$ (iii)] immediately. 
Hence it suffices to show that \textcircled{\scriptsize i} \& \textcircled{\scriptsize ii} $\R$ \textcircled{\scriptsize iii} $\R$ \textcircled{\scriptsize ii}. 
When $K=\mathbb Q(\zeta_9+\zeta_9^{-1})$, the assertion holds since $\mathbb Q(\zeta_9+\zeta_9^{-1})$ satisfies both of \textcircled{\scriptsize ii}, \textcircled{\scriptsize iii} 
by $\mathbb Q(\zeta_9+\zeta_9^{-1})=K_{0}$, $c_{\mathbb Q(\zeta_9+\zeta_9^{-1})}=\Delta_0=9$. 
We assume that $K\neq  \mathbb Q(\zeta_9+\zeta_9^{-1})$. \\[5pt]
[\textcircled{\scriptsize i} \& \textcircled{\scriptsize ii} $\R$ \textcircled{\scriptsize iii}]. 
Assume that \textcircled{\scriptsize i} holds. Then $\mathfrak C_K$ is principal by \cref{thm2}-(i).
There are three cases by \cref{expl}-[(i) $\LR$ (iii)]:
\begin{enumerate}
\item[(a)] $v_p(\Delta_t) \not \equiv 1 \bmod 3 \text{ for all } p$ and ($3 \nmid t \text{ or } t\equiv 12 \bmod 27$).
\item[(b)] $v_p(\Delta_t) \not \equiv 2 \bmod 3 \text{ for all }p$ and ($3 \nmid t \text{ or } t\equiv 12 \bmod 27$).
\item[(c)] $v_p(\Delta_t)\not \equiv 2 \bmod 3 \text{ for all } p\neq 3$ and ($t \equiv 0,6 \bmod 9$).
\end{enumerate}
Note that \textcircled{\scriptsize iii} states that only (b) or (c) holds.
On the other hand, \textcircled{\scriptsize ii} states that  
\begin{align} \label{pf1}
v_p(\Delta_t) \equiv v_p(c_K) \bmod 3 \text{ for all }p.
\end{align}
Since $K\neq  \mathbb Q(\zeta_9+\zeta_9^{-1})$, there exists $p \mid c_K,\neq 3$, which satisfies 
\begin{align} \label{pf2}
v_{p}(c_K)=1
\end{align}
by \cref{order}.
Hence (a) and (\ref{pf1}) do not hold simultaneously. \\[5pt]
[\textcircled{\scriptsize iii} $\R$ \textcircled{\scriptsize ii}].
Assume \textcircled{\scriptsize iii} holds. Then $\mathfrak C_K$ is principal since \textcircled{\scriptsize iii} is a part of the condition of \cref{expl}-(iii).
Therefore we have $\frac{\Delta_t}{c_K}$ or $\frac{\Delta_t^2}{c_K} \in \mathbb N^3$ by \cref{thm3}.
However, under \textcircled{\scriptsize iii}, $\frac{\Delta_t^2}{c_K} \in \mathbb N^3$ can not hold by (\ref{pf2}). Then the assertion is clear.

\end{document}

%% file: paper.tex
\usepackage{amsmath,amssymb,amsthm}
\usepackage{color}
\usepackage[final]{graphicx}
\usepackage{bm} 
\usepackage{ascmac} 
\usepackage[capitalize]{cleveref} 
\usepackage{comment} 
\usepackage{stmaryrd} 
\usepackage{ulem} 
\usepackage{url}
\usepackage[color,draft]{showkeys} 

\newtheorem{thm}{Theorem}[section]
\newtheorem{crl}[thm]{Corollary}

\newtheorem{prp}[thm]{Proposition}
\newtheorem{lmm}[thm]{Lemma}
\newtheorem{rmk}[thm]{Remark}

\newcommand{\R}{\Rightarrow}
\renewcommand{\L}{\Leftarrow}
\newcommand{\LR}{\Leftrightarrow}
\newcommand{\Llr}{\quad \Longleftrightarrow \quad}
\newcommand{\Lr}{\quad \Longrightarrow \quad}

\newcommand{\ra}{\rightarrow}
\newcommand{\mt}{\mapsto}

\makeatletter
\newcommand{\subjclass}[2][2010]{%
  \let\@oldtitle\@title%
  \gdef\@title{\@oldtitle\footnotetext{#1 \emph{Mathematics subject classification(s).} #2}}%
}
\newcommand{\keywords}[1]{%
  \let\@@oldtitle\@title%
  \gdef\@title{\@@oldtitle\footnotetext{\emph{Key words and phrases.} #1.}}%
}
\makeatother